\documentclass[12pt]{article}
\usepackage{amssymb, amsthm, amsmath}
\newtheorem{thm}{Theorem}[section]
\newtheorem{prop}{Proposition}[section]
\newtheorem{lem}{Lemma}[section]
\newtheorem{cor}{Corollary}[section]
\theoremstyle{definition}

\newtheorem{dfn}{Definition}[section]
\newtheorem{rem}{Remark}[section]
\newtheorem{ex}{Example}[section]

\def\D{{\overleftarrow{D}}}

\def\Sym{{\rm{Sym}}}
\def\AA{{\mathbb{A}}}

\def\ZZ{{\mathbb{Z}}}
\def\RR{{\mathbb{R}}}

\def\CC{{\mathbb{C}}}
\def\QQ{{\mathbb{Q}}}
\def\DD{{\mathbb{D}}}

\def\Gr{{\widehat{\rm Gr}}}
\def\Waf{{W_{\rm aff}}}
\def\Vaf{{V_{\rm aff}}}
\def\B{{\cal{B}}}
\def\Baf{{\cal{B}_{\rm aff}}}
\begin{document}
\title{Affine nil-Hecke algebras and braided differential structure 
on affine Weyl groups}
\date{}
\author{Anatol N. Kirillov and Toshiaki Maeno}
\maketitle
\begin{abstract}
We construct a model of the affine nil-Hecke algebra as a 
subalgebra of the Nichols-Woronowicz 
algebra associated to a Yetter-Drinfeld module over the 
affine Weyl group. We also discuss the Peterson isomorphism 
between the homology of the affine Grassmannian and the small 
quantum cohomology ring of the flag variety in terms of 
the braided differential calculus. 
\end{abstract}

\section*{Introduction}
The cohomology ring of the flag variety is a fundamental object of 
research in the study of the Schubert calculus. Fomin and the first author \cite{FK} 
gave a combinatorial 
model of the cohomology $H^*(Fl_n)$ 
ring of the flag variety of type $A$ as a commutative subalgebra 
of a quadratic algebra ${\cal E}_n$.  
It is remarkable that the algebra ${\cal E}_n$ has a natural quantum deformation 
${\cal E}_n^q$ so that ${\cal E}_n^q$ contains the quantum cohomology ring 
$QH^*(Fl_n)$ as a commutative subalgebra. 

It has been observed by Milinski and Schneider \cite{MS} and by Majid \cite{Maj2} 
that the defining relations of the Fomin-Kirillov 
quadratic algebra ${\cal E}_n$ are understandable from the viewpoint of a certain kind 
of braided Hopf algebra called the Nichols-Woronowicz algebra. 
Bazlov \cite{Ba} constructed 
the model of the coinvariant algebra of the 
finite Coxeter groups as a commutative subalgebra of the Nichols-Woronowicz algebra. 
At the same time, the nil-Coxeter algebra, which is dual to the coinvariant algebra, 
is also realized as a subalgebra of the Nichols-Woronowicz algebra. 

The braided analgue of the symmetric or exterior algebra was introduced by 
Woronowicz \cite{Wo} for the study of the differential forms on the quantum 
groups. For a given braided vector space $M$ over a field $K$ of characteristic zero, 
the braided analogue $\B(M)$ 
of the symmetric algebra of $M$ is defined to be the quotient of the 
free tensor algebra of $M$ by the kernel of the braided symmetrizer. 
It is known that the algebra $\B(M)$ is 
a braided graded Hopf algebra characterized by the following conditions: \\ 
(1) $\B^0(M)=K,$ \\ 
(2) $\B^1(M)=M=\{ \textrm{primitive elements in $\B(M)$} \},$ \\ 
(3) $\B^1(M)$ generates $\B(M)$ as an algebra. \\ 
The Hopf algebra characterized by the above conditions has been studied  
by Nichols \cite{Ni} and named the Nichols algebra 
by Andruskiewitsch and Schneider \cite{AS}. The study of the algebra 
$\B(M)$ from the viewpoint of the free braided differential calculus 
was developed by \cite{Maj1}. 
In this paper we will call $\B(M)$ 
the Nichols-Woronowicz algebra simply following \cite{Ba}. 

The aim of this paper is to construct the nil-Hecke algebra as a subalgebra 
of an extension of the Nichols-Woronowicz algebra $\Baf$ associated to a Yetter-Drinfeld module 
over the affine Weyl groups. Our construction is analogous to the one in \cite[Section 6]{Ba}. 

It is known that the affine Grassmannian $\Gr:=G(\CC((t)))/G(\CC[[t]])$ of a 
semisimple Lie group $G$ 
is homotopic to the loop group $\Omega K$ of the maximal 
compact subgroup $K\subset G.$ The homology $H_*(\Gr)\cong H_*(\Omega K)$ 
carries an associative algebra structure induced by the Pontryagin product. 
The strucuture of the Pontryagin ring $H_*(\Omega K)$ has been determined 
by Bott \cite{Bo}. The Schubert calculus for Kac-Moody flag varieties 
was studied by Kostant and Kumar \cite{KK} by using the nil-Hecke algebra. 
Peterson \cite{Pe} stated that the torus-equivariant homology $H^T_*(\Gr)$ 
of the affine Grassmannian is 
isomorphic to the so-called Peterson subalgebra of the affine nil-Hecke 
algebra. So our construction gives a model of $H^T_*(\Gr)$ as a commutative 
subalgebra of the Nichols-Woronowicz algebra $\Baf(S),$ see Theorem 3.1.  

Peterson \cite{Pe} also pointed out that the Pontryagin ring $H^T_*(\Gr)$ is 
isomorphic to the small quantum cohomology ring $QH^*_T(G/B)$ of the 
corresponding flag variety $G/B$ as an algebra after a suitable localization. 
The affine Bruhat operator acting on $H^T_*(\Gr)$ introduced by Lam and 
Shimozono \cite{LS} gives an explicit comparison between the multiplicative 
structure of $H^T_*(\Gr)$ and that of $QH^*_T(G/B).$ In this paper, 
we will realize the affine Bruhat operator as a braided differential operator 
acting on our algebra $\Baf.$ 

\subsection*{Acknowledgement.} 
The second author is supported by Grant-in-Aid for Scientific Research. 

\section{Affine nil-Hecke algebra}
Let $G$ be a simply-connected semisimple complex Lie group and $W$ its Weyl group. 
Denote by $\Delta$ the set of the roots. 
We fix the set $\Delta_+$ of the positive roots by choosing 
a set of simple roots $\alpha_1,\ldots,\alpha_r.$ 
The Weyl group $W$ acts on the weight lattice $P$ and the coroot lattice 
$Q^{\vee}$ of $G.$ 
The affine Weyl group $\Waf$ is generated by the affine reflections 
$s_{\alpha, k},$ $\alpha\in \Delta,$ $k\in \ZZ,$ with respect to the affine 
hyperplanes 
$H_{\alpha,k}:= \{ \lambda \in P\otimes \RR \; | \; \langle\lambda ,\alpha^{\vee}\rangle=k \}.$ 
The affine Weyl group is the semidirect product of $W$ and $Q^{\vee},$ i.e., 
$\Waf=W \ltimes Q^{\vee}.$ The affine Weyl group $\Waf$ is generated by the 
simple reflections $s_1:=s_{\alpha_1,0},\ldots,s_r:=s_{\alpha_r,0}$ and 
$s_0:=s_{\theta,1}$ 
where $\theta=-\alpha_0$ is the highest root. 
The affine Weyl group $W$ has the presentation as a Coxeter group as follows: 
\[ \Waf = \langle s_0,\ldots,s_r \; | \; s_0^2=\cdots =s_r^2=1, (s_is_j)^{m_{ij}}=1 
\rangle . \] 
\begin{dfn}
{\rm  
The affine nil-Coxeter algebra $\AA_0$ is the associative algebra generated by 
$\tau_0,\ldots,\tau_r$ subject to the relations 
\[ \tau_0^2=\cdots=\tau_r^2=0, \;\; 
(\tau_i\tau_j)^{[m_{ij}/2]}\tau_i^{\nu_{ij}}=(\tau_j\tau_i)^{[m_{ij}/2]}\tau_j^{\nu_{ij}} , \] 
where $\nu_{ij}:=m_{ij}-2[m_{ij}/2].$}
\end{dfn}

For a reduced expression $x=s_{i_1}\cdots s_{i_l}$ of an element $x\in \Waf,$ 
the element $\tau_x:=\tau_{i_1}\cdots \tau_{i_l}\in \AA_0$ is independent 
of the choice of the reduced expression of $x.$ It is known that $\{ \tau_x \}_{x\in \Waf}$ 
form a linear basis of $\AA_0.$ 

The nil-Coxeter algebra $\AA_0$ acts on $S:=\Sym P_{\QQ}$ via 
\[ \tau_0(f):=\partial_{\alpha_0}(f)=-(f-s_{\theta,0}f) / \theta , \] 
\[ \tau_i(f):= \partial_{\alpha_i}(f)=(f-s_{\alpha_i,0}f)/\alpha_i, \; \; i=1,\ldots, r, \] 
for $f\in S.$ 
\begin{dfn}
{\rm (\cite{KK}) The nil-Hecke algebra $\AA$ is defined to be the cross product 
$\AA_0 \ltimes S,$ where the cross relation is given by 
\[ \tau_i f= \partial_{\alpha_i}(f)+s_i(f)\tau_i \;\; f\in S, i=1,\ldots, r. \] }
\end{dfn}

Here, we summarize some known results on the homology of 
the affine Grassmannian. 
The affine Grassmannian $\Gr:=G(\CC((t)))/G(\CC[[t]])$ 
is homotopic to the loop group $\Omega K$ of the maximal 
compact subgroup $K\subset G.$ 
Let $T\subset G$ be the maximal torus. 
An associative algebra structure on the 
$T$-equivariant homology group 
$H^T_*(\Gr)\cong H^T_*(\Omega K)$ is induced from the group multiplication 
\[ \Omega K \times \Omega K \rightarrow \Omega K. \] 
It is known that the algebra $H^T_*(\Gr)$ is commutative. 
The algebra $H^T_*(\Omega K)$ is called the Pontryagin ring. 

We regard the $T$-equivariant homology $H_*^T(\Gr)$ as an $S$-algebra 
by identifying $S=H_T^*(pt).$ The diagonal embedding 
\[ \Omega K \rightarrow \Omega K \times \Omega K \] 
induces a coproduct on $H_*^T(\Gr).$ 
\begin{prop}
{\rm (\cite{Pe})} The $T$-equivariant homology $H_*^T(\Gr)$ is isomorphic to 
the centralizer $Z_{\AA}(S)$ of $S$ in $\AA$ as Hopf algebras. 
\end{prop} 

\section{Nichols-Woronowicz algebra for affine Weyl groups} 
We briefly recall the construction of the Nichols-Woronowicz algebra 
associated to a braided vector space. Let $M$ be a vector space over a field of 
characteristic zero and 
$\psi:M^{\otimes 2} \rightarrow M^{\otimes 2}$ be a fixed linear endomorphism satisfying 
the braid relations $\psi_{i}\psi_{i+1}\psi_{i}=\psi_{i+1}\psi_{i}\psi_{i+1}$ where 
$\psi_i:M^{\otimes n} \rightarrow M^{\otimes n}$ is a linear endomorphism 
obtained by applying $\psi$ to the $i$-th and $(i+1)$-st components. 
Denote by $s_i$ the simple transposition $(i,i+1)\in S_n.$ 
For any reduced expression
$w=s_{i_1}\cdots s_{i_l}\in S_n,$ the endomorphism 
$\Psi_w=\psi_{i_1}\cdots \psi_{i_l}:M^{\otimes n} 
\rightarrow M^{\otimes n}$ is well-defined. 
The Woronowicz symmetrizer \cite{Wo} is given by 
$\sigma_n := \sum_{w\in S_n} \Psi_w .$ 

\begin{dfn} {\rm (\cite{Wo})} 
The Nichols-Woronowicz algebra associated to a braided vector space $M$ is 
defined by 
\[ \B(M):= \bigoplus_{n\geq 0} M^{\otimes n}/{\rm Ker}(\sigma_n), \] 
where $\sigma_n : M^{\otimes n} \rightarrow M^{\otimes n}$ is the Woronowicz 
symmetrizer. 
\end{dfn}

\begin{dfn}
A vector space $M$ is called a Yetter-Drinfeld module over a group $\Gamma,$ if 
the following conditions are satisfied: \\ 
$(1)$ $M$ is a $\Gamma$-module, \\
$(2)$ $M$ is $\Gamma$-graded, i.e. $M=\bigoplus_{g\in \Gamma}M_g,$ where $M_g$ is a 
linear subspace of $M,$ \\ 
$(3)$ for $h\in \Gamma$ and $v\in M_g,$ $h(v)\in M_{hgh^{-1}}.$ 
\end{dfn} 

The Yetter-Drinfeld module $M$ over a group $\Gamma$ is naturally braided 
with the braiding $\psi:M^{\otimes 2} \rightarrow M^{\otimes 2}$ defined by 
$\psi(a \otimes b)= g(b) \otimes a$ for $a \in M_g$ and $b\in M.$ 

In the following we are interested in the Yetter-Drinfeld module over the 
affine Weyl groups $\Waf.$ 
Denote by $t_{\lambda} \in \Waf$ the translation by $\lambda \in Q^{\vee}.$ 
We define a Yetter-Drinfeld module $\Vaf$ over $\Waf$ by 
\[ \Vaf := \bigoplus_{\alpha\in \Delta, k\in \ZZ} \QQ \cdot [\alpha,k]/
( [\alpha,k]+[-\alpha,-k]) , \] 
where the $\Waf$ acts on $\Vaf$ by 
\[ w [\alpha,k] := [w(\alpha),k], \;\; w\in W, \;\;\; 
t_{\lambda}[\alpha,k] := [\alpha,k+(\alpha,\lambda)] , \;\; \lambda \in Q^{\vee}. \] 
The $\Waf$-grading is given by $\deg_{\Waf}([\alpha,k]):=s_{\alpha,k}.$ 
Then it is easy to check the conditions in Definition 2.1. 
Now we have the Nichols-Woronowicz algebra $\Baf:=\B(\Vaf)$ associated to the Yetter-Drinfeld 
module $\Vaf.$ 

Let $\B_W$ be the Nichols-Woronowicz algebra associated to 
the Yetter-Drinfeld module $V=\oplus_{\alpha \in \Delta} \QQ \cdot [\alpha]/
([\alpha]+[-\alpha])$ as in \cite[Section 4]{Ba}. 
\begin{lem} 
$(1)$ We have a surjective homomorphism $\pi:\Baf \rightarrow \B_W,$ 
$\pi([ \alpha  ,k ]):=  [ \alpha ] .$ \\ 
$(2)$ The algebra $\Baf$ acts on $S$ via $[\alpha,k] f= \partial_{\alpha}(f)$ for all 
$k\in \ZZ.$ 
\end{lem}
{\it Proof.} 
$(1)$ Denote by $\psi$ and $\bar{\psi}$ the braidings on $\Vaf$ and $V$ respectively.  
Let $\tilde{\pi}:\oplus_n \Vaf^{\otimes n} 
\rightarrow \oplus_n V^{\otimes n}$ be the lift of $\pi.$
Since 
\[ \psi([\alpha,k] \otimes [\beta,l])=
[s_{\alpha}(\beta), l-\langle \alpha^{\vee},\beta \rangle k] \otimes [ \alpha, k ] \] 
and 
$\bar{\psi}([\alpha] \otimes [\beta])=
[s_{\alpha}(\beta)]\otimes [\alpha],$ the map $\tilde{\pi}$ 
sends the kernel of the braided symmetrizer $\sigma_n$ of $\Vaf^{\otimes n}$ to that of $V^{\otimes n}.$ \\ 
$(2)$ In \cite{Ba}, it is shown that the algebra $\B_W$ acts on the coinvariant 
algebra $S_W$ via $[\alpha] \mapsto \partial_{\alpha}.$ Let $S^W$ be the $W$-invariant 
subalgebra of $S.$ Then we have the decomposition $S=S^W \otimes S_W.$ The operator $\partial_{\alpha}$ 
extends $S^W$-linearly to the operator on $S.$ Hence $\B_W$ acts on $S.$ 
We have seen the existence of the natural projection $\pi$ from $\Baf$ to $\B,$ 
so $\pi$ induces the action of $\Baf$ on $S.$ 
\bigskip \\ 
Let us define the extension $\Baf (S)=\Baf \ltimes S$ by the cross 
relation 
\[ [\alpha,k] f = \partial_{\alpha}f +s_{\alpha,0}(f)[\alpha,k] , \;\; 
[\alpha,k] \in \Vaf, f\in S. \] 

\begin{prop}
There exists a homomorphism $\varphi : \AA \rightarrow \Baf (S)$ given by 
$\tau_0 \mapsto [\alpha_0,-1],$ $\tau_i \mapsto [\alpha_i,0],$ $i=1,\ldots,r,$ and 
$f \mapsto f,$ $f\in S.$  
\end{prop} 
{\it Proof.} It is enough to check the Coxeter relations among 
$\varphi(\tau_0),\ldots ,\varphi(\tau_r)$ in $\Baf(S)$ 
based on the classification of the affine root systems. 
This is done by the direct computation of the symmetrizer for the subsystems of 
rank 2 in the similar manner to \cite[Section 6]{Ba}. 

\begin{ex} {\rm 
Here we list the Coxeter relations in $\Baf$ involving $[\theta,1]=-[\alpha_0,-1]$ 
for the root systems of rank 2. Let $(\varepsilon_1,\ldots,\varepsilon_r)$ be an orthonormal 
basis of the $r$-dimensional Euclidean space. Put 
$[ij,k]:=[\varepsilon_i-\varepsilon_j,k],$ 
$[\overline{ij},k]:=[\varepsilon_i+\varepsilon_j,k],$ 
$[i,k]:=[\varepsilon_i,k]$ and $[\alpha]:=[\alpha,0].$ \\ 
(i) (Type $A_2$ case) 
\[ [13,1][23][13,1]+[23][13,1][23]=0, \;\; [13,1][12][13,1]+[12][13,1][12]=0 \] 
(ii) (Type $B_2$ case)
\[ [\overline{12},1][2][\overline{12},1][2]=[2][\overline{12},1][2][\overline{12},1] \]
(iii) (Type $G_2$ case) Let $\alpha_1,\alpha_2$ be the simple roots for 
$G_2$-system. We assume that $\alpha_1$ is a short root and $\alpha_2$ is a long one. 
Then we have $\theta=3\alpha_1+2\alpha_2.$ 
\[ [\theta,1][\alpha_2][\theta,1]+[\alpha_2][\theta,1][\alpha_2]=0. \] 
} 
\end{ex}

\section{Model of nil-Hecke algebra} 
The connected components of $P\otimes \RR \setminus \cup_{\alpha\in \Delta_+,k\in \ZZ}H_{\alpha ,k}$ 
are called alcoves. The affine Weyl group $\Waf$ acts on the set of the alcoves 
simply and transitively. 
\begin{dfn}{\rm (\cite{LP})} $(1)$ A sequence $(A_0,\ldots,A_l)$ of alcoves $A_i$ is called 
an alcove path if $A_i$ and $A_{i+1}$ have a common wall and $A_i \not= A_{i+1}.$ \\ 
$(2)$ An alcove path $(A_0,\ldots,A_l)$ is called reduced if the length $l$ 
of the path is minimal among all alcove paths connecting $A_0$ and $A_l.$ \\ 
$(3)$ We use the symbol $A_i \stackrel{\beta,k}{\longrightarrow} A_{i+1}$ when 
$A_i$ and $A_{i+1}$ have a common wall of the form $H_{\beta,k}$ and the direction 
of the root $\beta$ is from $A_i$ to $A_{i+1}.$ 
\end{dfn} 

The alcove $A^{\circ}$ defined by the inequalities 
$\langle \lambda,\alpha_0^{\vee}\rangle \geq -1$ and 
$\langle \lambda,\alpha_i^{\vee}\rangle \geq 0,$ $i=1,\ldots,r,$ is 
called the fundamental alcove. 
For a reduced alcove path $\gamma:A_0=A^{\circ} \stackrel{\beta_1,k_1}{\longrightarrow} \cdots 
\stackrel{\beta_l,k_l}{\longrightarrow} A_l,$ we define an element 
$[\gamma]\in \Baf$ by 
\[ [ \gamma ] := [-\beta_1,-k_1] \cdots [-\beta_l,-k_l] . \] 
When $A_l=x^{-1}(A^{\circ})$ for $x\in \Waf,$ we will also use the symbol $[x]$ instead of $[\gamma],$ 
since $[\gamma]$ depends only on $x$ thanks to the Yang-Baxter relation. 

For a braided vector space $M,$ 
it is known that an element $a\in M$ acts on $\B(M^*)$ as a braided 
differential operator (see \cite{Ba}, \cite{Maj1}). 
Let us identify $M^*$ with $M$ via the $\Waf$-invariant 
inner product $(\; , \; )$ given by 
\[ ([\alpha,k],[\beta ,l])= \left\{ 
\begin{array}{cc}
1, & \textrm{if $\alpha=\beta$ and $k=l,$} \\
0, & \textrm{otherwise,}  
\end{array}
\right. \]
for $\alpha,\beta \in \Delta_+,$ $k,l\in \ZZ.$ 
In our case, the differential operator $\D_{[\alpha,k]},$ 
$[\alpha,k]\in \Vaf,$ acting 
from the right is determined by the following characterization: \\ 
(0) $(c)\D_{[\alpha,k]}=0,$ $c\in \QQ,$ \\ 
(1) $([\alpha,k])\D_{[\beta ,l]} = ([\alpha,k],[\beta ,l]),$ \\ 
(2) $(FG)\D_{[\alpha,k]}=F(G\D_{[\alpha,k]})+(F\D_{[\alpha,k]})s_{\alpha,k}(G),$ \\ 
for $\alpha,\beta \in \Delta,$ $k,l\in \ZZ,$ $F,G\in \Baf.$ 
The operator $\D_{[\alpha,k]}$ extends to the one acting on $\Baf(S)$ by the 
commutation relation $f\cdot\D_{[\alpha,k]}=\D_{[\alpha,k]}\cdot s_{\alpha,k}(f),$ 
$f\in S.$ 

We use the abbreviation $\D_0:=\D_{[\alpha_0,-1]},$ $\D_i:=\D_{[\alpha_i,0]},$ 
$i=1,\ldots,r.$ 
For $x\in \Waf,$ fix a reduced decomposition $x=s_{i_1}\cdots s_{i_l}.$ 
We define the corresponding 
braided differential operator $\D_x$ acting on $\Baf$ by the formula 
\[ \D_x := \D_{i_l}\cdots \D_{i_1} , \] 
which is also independent of the choice of the reduced decomposition of $x$ 
because of the braid relations. 
\begin{lem}
For $x\in \Waf,$ take a reduced alcove path $\gamma$ from the fundamental 
alcove $A^{\circ}$ to $x^{-1}(A^{\circ}).$ Then, we have $([\gamma])\D_x=1.$ 
\end{lem}
{\it Proof.} Let us take a reduced path 
\[ \gamma: A_0=A^{\circ} \stackrel{\beta_1,k_1}{\longrightarrow} A_1 
\stackrel{\beta_2,k_2}{\longrightarrow} \cdots 
\stackrel{\beta_l,k_l}{\longrightarrow} A_l=x^{-1}(A^{\circ}). \] 
Define a sequence $\sigma_1,\ldots,\sigma_l \in \Waf$ inductively by 
\[ \sigma_1:=s_{\beta_1,k_1}, \; 
\sigma_{j+1}:= \sigma_j s_{\beta_{j+1},k_{j+1}} \sigma_j. \] 
Then it is easy to see that $\sigma_{\nu}(A_j)\not= A^{\circ},$ $1\leq \nu \leq j-1,$ 
$\sigma_j(A_j)=A^{\circ}$ and the walls 
$\sigma_j(H_{\beta_{j+1},k_{j+1}})$ are corresponding to simple roots. 
Hence, $\sigma_1,\ldots,\sigma_l$ are 
simple reflections. This sequence gives a reduced expression 
$x=\sigma_l \cdots \sigma_1.$ 
Put $\sigma_i=s_{\alpha_{i_j}}.$ 
Since the direction of $\beta_{j+1}$ is chosen to be 
from $A_j$ to $A_{j+1},$ we have 
\[ [\gamma] \D_x= ([\beta_1,k_1])\D_{i_1}\cdot (\sigma_1([\beta_2,k_2]))\D_{i_2} 
\cdots (\sigma_{l-1}([\beta_l,k_l]))\D_{i_l} =1. \] 

\begin{ex}
($A_2$-case) 
The standard realization is given by $\alpha_1=\varepsilon_1-\varepsilon_2,$ 
$\alpha_2=\varepsilon_2-\varepsilon_3,$ 
$\alpha_0=\varepsilon_3-\varepsilon_1.$ 
Consider the translation $t_{\alpha_1}$ by the simple root $\alpha_1.$ 
If we take a reduced path 
\[ \gamma: A_0=A^{\circ} \stackrel{-\alpha_2,0}{\longrightarrow} A_1 
\stackrel{\alpha_1,1}{\longrightarrow} A_2 
\stackrel{-\alpha_0,1}{\longrightarrow} A_3 
\stackrel{\alpha_1,2}{\longrightarrow} A_4=t_{\alpha_1}(A^{\circ}) , \] 
then we have $[\gamma]=[23][21,-1][31,-1][21,-2].$ 
On the other hand, the differential operator corresponding to $t_{-\alpha_1}$ 
is given by $\D_2\D_0\D_2\D_1,$ where $\D_0=\D_{[31,-1]},\D_1=\D_{[12]},
\D_2=\D_{[23]}.$ 
It is easy to check by direct computation 
\[ ([23][21,-1][31,-1][12,2])\D_2\D_0\D_2\D_1=1 . \] 
\end{ex}

\begin{thm}
The algebra homomorphism $\varphi: \AA \rightarrow \Baf (S)$ is injective. 
\end{thm}
{\it Proof.} The nil-Hecke algebra $\AA$ is also $\Waf$-graded. 
Since the homomorphism $\varphi:\AA \rightarrow \Baf (S)$ preserves the $\Waf$-grading, 
it is enough to check $\varphi(\tau_x) \not=0,$ for $x\in \Waf$ in order to show the injectivity 
of $\varphi.$ On the other hand, $\Baf^{op}$ acts on $\Baf$ itself via the 
braded differential operators. 
Let $\gamma$ be a reduced alcove path from $A^{\circ}$ to $x^{-1}(A^{\circ}).$ 
Then we have $([\gamma])\D_x=1$ from Lemma 3.1. This shows $\D_x\not=0,$ so 
$\varphi(\tau_x)\not=0.$ 
\medskip 

This theorem implies the following (see Proposition 1.1): 
\begin{cor}
The $T$-equivariant Pontryagin ring $H_*^T(\Gr)$ is a subalgebra 
of $\Baf(S).$ 
\end{cor} 
By taking the non-equivariant limit, we also have: 
\begin{cor}
The Pontryagin ring $H_*(\Gr)$ is a subalgebra of $\Baf.$ 
\end{cor}

\section{Affine Bruhat operators} 
We denote by $x \rightarrow y$ the cover relation in the Bruhat 
ordering of $\Waf,$ i.e. $y=xs_{\alpha,k}$ for some $\alpha \in \Delta$ 
and $k\in \ZZ,$ and $l(y)=l(x)+1.$ 

We will use some terminology from \cite{LS}. Denote by $\tilde{Q}$ the set of 
antidominant elements in $Q^{\vee}.$ An element $x\in \Waf$ can be 
expressed uniquely as a product of form $x=wt_{v\lambda}\in \Waf$ 
with $v,w \in W,$ $\lambda \in \tilde{Q}.$ 
We say that $x=wt_{v\lambda}$ belongs to the "$v$-chamber". 
An element $\lambda \in \tilde{Q}$ is called superregular when 
$|\langle \lambda, \alpha \rangle | > 2(\# W) +2$ for all $\alpha \in \Delta_+.$ 
If $\lambda\in \tilde{Q}$ is superregular, then $x=wt_{v\lambda}$ 
is called superregular. 
The subset of superregular elements in $\Waf$ is denoted by $\Waf^{\rm sreg}.$ 
We say that a property holds for sufficiently superregular elements 
$\Waf^{\rm ssreg}\subset \Waf$ if there is a positive constant $k\in \ZZ$ such that 
the property holds for all $x\in \Waf^{\rm sreg}$ satisfying the following condition: 
\[ y \in \Waf,\; y<x,\; \textrm{and} \; l(x)-l(y)<k \Rightarrow y\in \Waf^{\rm sreg}. \] 
The meaning of $\Waf^{\rm ssreg}$ depends on the context, see \cite[Section 4]{LS} 
for the details. 
For $v\in W,$ consider the $S$-submodule $M^{\rm ssreg}_v$ in $\Baf$ generated by 
the sufficiently superregular elements $[x]$ where $x$ belongs to the $v$-chamber. 

\begin{lem} Let $x\in \Waf.$ 
For $\alpha \in \Delta$ and $k\in \ZZ_{> 0},$ we have 
\[ [x] \D_{[\alpha,k]} = \left\{ 
\begin{array}{cc}
[xs_{\alpha,k}] , & \textrm{if $l(x)=l(xs_{\alpha,k})+1,$} \\ 
0, & \textrm{otherwise.}
\end{array}
\right. \] 
\end{lem}
{\it Proof.} The fundamental 
alcove $A^{\circ}$ is contained in the region $\{ \lambda \in P\otimes \RR | 
\langle \lambda,\alpha^{\vee} \rangle < k \}$ for $\alpha \in \Delta$ and 
$k\in \ZZ_{> 0}.$ 
Let us choose any reduced path $\gamma:A_0 \stackrel{\beta_1,k_1}{\longrightarrow} \cdots 
\stackrel{\beta_l,k_l}{\longrightarrow} A_l=x^{-1}(A^{\circ})$ with $k_i \geq 0.$ 
If $l(x)>l(xs_{\alpha,k}),$ then 
$(\beta_i,k_i)=(\alpha,k)$ for some $i.$ Take the largest $i$ and consider the path 
\[ \gamma':A_0 \stackrel{\beta_1,k_1}{\longrightarrow} \cdots 
\stackrel{\beta_{i-1},k_{i-1}}{\longrightarrow} A_{i-1} 
\stackrel{\beta'_{i+1},k'_{i+1}}{\longrightarrow} s_{\alpha,k}(A_{i+1}) 
\stackrel{\beta'_{i+2},k'_{i+2}}{\longrightarrow} \cdots \] 
\[ \cdots 
\stackrel{\beta'_l,k'_l}{\longrightarrow} s_{\alpha,k}(A_l)=
s_{\alpha,k}x^{-1}(A^{\circ})=(xs_{\alpha,k})^{-1}(A^{\circ}), \] 
where $(\beta'_j,k'_j)$ is determined by the condition 
$s_{\alpha,k}(H_{\beta_j,k_j})=H_{\beta'_j,k'_j}.$ 
If $l(x)=l(xs_{\alpha,k})+1,$ then the path $\gamma'$ 
is a reduced path. 
In this case, we have 
$[x] \D_{[\alpha,k]}=[xs_{\alpha,k}].$ If $l(x) > l(xs_{\alpha,k})+1,$ 
the above path $\gamma'$ is not reduced and $[x] \D_{[\alpha,k]}=0.$ 
When $l(x)<l(xs_{\alpha,k}),$ the element $[\alpha,k]$ does not appear 
in the monomial $[\gamma],$ so we have $[x] \D_{[\alpha,k]}=0.$ 

\begin{prop} {\rm (\cite[Proposition 4.1]{LS})} Let $\lambda \in \tilde{Q}$ be superregular. 
For $x=wt_{v\lambda}$ and $y=xs_{v\alpha,-n}$ with $v,w\in W,$ 
we have the cover relation $y \rightarrow x$ if and only if one of the following conditions 
holds: \\ 
$(1)$ $l(wv)=l(wvs_{\alpha})-1$ and $n=\langle \lambda, \alpha \rangle,$ 
giving $y=ws_{v(\alpha)}t_{v(\lambda)},$ \\ 
$(2)$ $l(wv)=l(wvs_{\alpha})+\langle \alpha^{\vee},2\rho \rangle -1$ and 
$n=\langle \lambda, \alpha \rangle +1,$ 
giving $y=ws_{v(\alpha)}t_{v(\lambda+\alpha^{\vee})},$ \\ 
$(3)$ $l(v)=l(vs_{\alpha})+1$ and $n=0,$ giving 
$y=ws_{v(\alpha)}t_{vs_{\alpha}(\lambda)},$ \\ 
$(4)$ $l(v)=l(vs_{\alpha})-\langle \alpha^{\vee},2\rho \rangle +1$ and $n=-1,$ 
giving $y=ws_{v(\alpha)}t_{vs_{\alpha}(\lambda+\alpha^{\vee})}.$
\end{prop}
In \cite{LS}, the first kind of the conditions (1) and (2) are called the near relation 
because 
$x$ and $y$ belong to the same chamber. In this paper we denote the near 
relation by $y \rightarrow_{near} x.$ 

The affine Bruhat operator $B^{\mu}:S\langle \Waf^{\rm ssreg} \rangle 
\rightarrow S\langle \Waf^{\rm sreg}\rangle,$ $\mu \in P,$ due to Lam and Shimozono \cite[Section 5]{LS} 
is an $S$-linear map defined by the formula 
\[ B^{\mu}(x) =(\mu - wv \mu) x + 
\sum_{\alpha\in \Delta_+}\sum_{xs_{v(\alpha),k}\rightarrow_{near} x}
\langle \alpha^{\vee}, \mu \rangle  xs_{v(\alpha),k}  \] 
for $x=wt_{v\lambda}\in \Waf^{\rm ssreg}.$ 
We also introduce the operator $\beta^{\mu}_v,$ $\mu \in P,$ acting on each $M^{\rm ssreg}_v$ by 
\[ \beta^{\mu}_v([x]):=(\mu - wv \mu) [x] + 
[x] \sum_{\alpha\in \Delta_+,k>1} 
\langle \alpha^{\vee}, \mu \rangle \D_{[v(\alpha),k]}, \] 
where $x=wt_{v\lambda} \in \Waf^{\rm ssreg}.$ 
Denote by $\Waf^{\rm ssreg}(v)$ the subset of $\Waf$ consisting of the 
superregular elements belonging to the $v$-chamber. 
Fix a left $S$-module 
isomorphism 
\[ \begin{array}{cccc} 
\iota : & S \langle \Waf^{\rm ssreg}(v) \rangle & \rightarrow & M_v^{\rm ssreg} \\ 
 & x & \mapsto & [x] . 
\end{array} 
\] 
\begin{prop}
For each $v\in W$ and a sufficiently superregular element $x \in \Waf^{\rm ssreg}(v),$ 
\[ \beta^{\mu}_v([x])=\iota (B^{\mu}(x)). \] 
\end{prop}
{\it Proof.} This can be shown by using Lemma 4.1 and Proposition 4.1. 
\[ \beta^{\mu}_v( [x])=(\mu - wv \mu) [x] + 
[x] \sum_{\alpha\in \Delta_+,k>1} 
\langle \alpha^{\vee}, \mu \rangle \D_{[v(\alpha),k]} \] 
\[ =(\mu - wv \mu) [x] + 
\sum_{\alpha\in \Delta_+}\sum_{k>1,l(xs_{[v(\alpha),k}])=l(x)-1}
\langle \alpha^{\vee}, \mu \rangle [xs_{v(\alpha),k}] \]
\[ =(\mu - wv \mu) [x] + 
\sum_{\alpha\in \Delta_+}\sum_{xs_{v(\alpha),k}\rightarrow_{near} x}
\langle \alpha^{\vee}, \mu \rangle  [xs_{v(\alpha),k}]  = 
 \iota (B^{\mu}(x)). \]
\begin{rem}
{\rm In \cite{KM} the authors introduced the quantization operators $\eta_{\alpha}$ acting 
on the model of $H^*(G/B)\otimes \CC[q_1,\ldots,q_r]$ realized as a subalgebra 
of $\B_W \otimes \CC[q_1,\ldots,q_{n-1}].$ For a superregular element $\lambda \in \tilde{Q}$ 
and $w\in W,$ consider a homomorphism $\theta_w^{\lambda}$ from the $\lambda$-small 
elements (see \cite[Section 5]{LS}) of $H^*(G/B)\otimes \CC[q]$ to $\Baf$ 
defined by 
\[ \theta_w^{\lambda}(q^{\mu}\sigma^v):= [vw^{-1}t_{w(\lambda+\mu)}], \] 
where $\sigma^v$ is the Schubert class of $G/B$ corresponding to $v\in W$ and 
$q^{\mu}=q_1^{\mu_1}\cdots q_r^{\mu_r}$ for $\mu=\sum_{i=1}^r\mu_i\alpha_i^{\vee}.$ 
The following is an interpretation of the formula of \cite[Proposition 5.1]{LS} 
in our setting:  
\[ \theta_w^{\lambda}(\eta_{\alpha}(\sigma))=
\beta_w^{\varpi_{\alpha}}(\theta_w^{\lambda}(\sigma)). \] 
}
\end{rem}

\section{Quadratic relations}
For $\alpha\in \Delta_+$ and $v\in W,$ let us define the operator $\DD_v(\alpha)$ by 
\[ \DD_v(\alpha):=\sum_{k>1}  \D_{[v(\alpha),k]}. \] 
Then we have 
\[ \beta^{\mu}_v([x])= (\mu - wv \mu) [x] +
[x]\sum_{\alpha \in \Delta_+} \langle \alpha^{\vee}, \mu \rangle \DD_v(\alpha) . \] 
In the following, we discuss the relations among the operators $\DD_v(\alpha),$ $\alpha\in \Delta_+,$ 
for the root system of type $A_{n-1}.$ For simplicity, we consider only non-equivariant case with $v=id.$ 
Take the standard realization of the $A_{n-1}$-system: 
\[ \Delta= \{ \varepsilon_i-\varepsilon_j \; | \; 1\leq i,j \leq n, i\not=j \} . \] 
Put $\DD(ij):=\DD_{id.}(\varepsilon_i-\varepsilon_j)$ for $1\leq i<j \leq n,$ and 
$\DD(ij):=-\DD(ji)$ for $i>j.$ 
In this situation, we have a formula for the non-equivariant limit 
$\bar{\beta}^{\varepsilon_i}_{id.}$ of the operator $\beta^{\varepsilon_i}_{id.}:$ 
\[ \bar{\beta}^{\varepsilon_i}_{id.}=\sum_{j\not=i}\DD(ij) . \] 
Note that this formula is analogous to the definition of the Dunkl elements in \cite{FK}. 

Let $T_i,$ $1\leq i \leq n-1,$ be linear operators on $M^{\rm ssreg}$ defined 
by $T_i([x]):=[xt_{\alpha_i}],$ where $x\in \Waf$ and 
$\alpha_i=\varepsilon_i-\varepsilon_{i+1}.$ 
It is easy to check from Proposition 4.1 that 
$(T_i[x])\DD(jk)=T_i([x]\DD(jk)).$ 
Our next goal is to show that the operators $\DD(ij)$ satisfy 
the defining relations of the quantum deformation ${\cal E}_n^q$ of the 
Fomin-Kirillov quadratic algebra \cite{FK}. 
\begin{prop} 
{\rm (i)} For $1\leq i<j \leq n,$ we have 
\[ \DD(ij)^2= \left\{ \begin{array}{cc}
T_i, & \textrm{if $j=i+1,$} \\ 
0, & \textrm{otherwise.} 
\end{array} \right. 
\]
{\rm (ii)} If $\{i,j\} \cap \{ k,l \}=\emptyset,$ then we have $\DD(ij)\DD(kl)=\DD(kl)\DD(ij).$ \\ 
{\rm (iii)} For $1\leq i,j \leq n,$ $i\not= j,$ we have 
\[ \DD(ij)\DD(jk)+\DD(jk)\DD(kl)+\DD(ki)\DD(ij)=0. \] 
\end{prop}
{\it Proof.} First of all, let us check the equality (i). We have 
\[ \DD(ij)^2= \sum_{k,l>1}\D_{[ij,k]} \D_{[ij,l]} . \] 
Let $\lambda \in \tilde{Q}$ be sufficiently superregular. 
For $x=wt_{\lambda}\in \Waf,$ assume that $[x]\D_{[ij,k]} \D_{[ij,l]} \not=0.$ 
Then we have the arrows 
$xs_{ij,k} \rightarrow_{near} x$ 
and $xs_{ij,k}s_{ij,l} \rightarrow_{near} xs_{ij,k}$ 
in the Bruhat ordering. From the conditions (1) and (2) in 
Proposition 4.1, one of the following conditions holds: \\ 
Case (1): $k=-\langle \lambda ,\varepsilon_i-\varepsilon_j \rangle$ and 
$l(w)=l(ws_{ij})-1,$ \\ 
Case (2): $k=-\langle \lambda, \varepsilon_i-\varepsilon_j \rangle -1$ and 
$l(w)=l(ws_{ij})+\langle \varepsilon_i-\varepsilon_j,2\rho \rangle -1.$ \\ 
In Case (1), since the arrow 
$xs_{ij,k}s_{ij,l}=ws_{ij}t_{\lambda}s_{ij,l} \rightarrow_{near} xs_{ij,k}$ 
must come from the condition (2) of Proposition 4.1, we have $\langle \varepsilon_i-\varepsilon_j,2\rho \rangle -1=1.$ 
This equality implies that $\varepsilon_i-\varepsilon_j$ is a simple root $\alpha_i,$ and 
we get 
\[ [x]\DD(i \; i+1)^2 =
[x]\D_{[\alpha_i,-\langle \lambda ,\alpha_i \rangle]}\D_{[\alpha_i,-\langle \lambda ,\alpha_i \rangle-1]}
=[xt_{\alpha_i}]=T_i[x]. \] 
In Case (2), since 
the arrow $xs_{ij,k}s_{ij,l}=
ws_{ij}t_{\lambda+\varepsilon_i-\varepsilon_j}s_{ij,l} \rightarrow_{near} xs_{ij,k}$ 
comes from the condition (1) of Proposition 4.1, we again obtain 
$\langle \varepsilon_i-\varepsilon_j,2\rho \rangle -1=1$ and 
$\varepsilon_i-\varepsilon_j=\alpha_i.$ Hence we get 
\[ [x]\DD(i \; i+1)^2 = 
[x]\D_{[\alpha_i,-\langle \lambda ,\alpha_i \rangle-1]}\D_{[\alpha_i,-\langle \lambda ,\alpha_i \rangle-2]}
=[xt_{\alpha_i}]=T_i[x]. \] 
If $j\not=i+1,$ we have 
$\DD(ij)^2=0.$ 
The relations (ii) and (iii) follow from the identities  
$[ij,a][kl,b]=[kl,b][ij,a]$ for $\{ i,j\} \cap \{ k,l \} =\emptyset,$ and 
\[ [ij,a][jk,b]+[jk,b][ki,-a-b]+[ki,-a-b][ij,a]=0 \] 
in $\Baf.$ 

\begin{rem}
The operators $\DD_v(\alpha)$ induce the quantum Bruhat representation of ${\cal E}_n^q$ 
via $\theta_v^{\lambda}.$ 
\end{rem}

Research Institute for Mathematical Sciences \\ 
Kyoto University \\
Sakyo-ku, Kyoto 606-8502, Japan \\
e-mail: {\tt kirillov@kurims.kyoto-u.ac.jp} \\
URL: {\tt http://www.kurims.kyoto-u.ac.jp/\textasciitilde kirillov} 
\bigskip \\ 
Department of Electrical Engineering, \\
Kyoto University, \\ 
Sakyo-ku, Kyoto 606-8501, Japan \\ 
e-mail: {\tt maeno@kuee.kyoto-u.ac.jp}
\end{document}